\documentclass[12pt]{article}

\usepackage{graphicx}

\title{\bf Normally Regular Digraphs}
\author{
Leif K. J\o rgensen \\
\small Dept. of Mathematical Sciences\\
\small Aalborg University\\
\small Fr. Bajers Vej 7,
\small 9220 Aalborg, Denmark.\\
\small \tt leif@math.aau.dk}

\usepackage{amsthm,amsmath,amssymb}
\newtheorem{thm}{Theorem}
\newtheorem{conj}{Conjecture}
\newtheorem{lemma}[thm]{Lemma}
\newtheorem{prop}[thm]{Proposition}
\newtheorem{cor}[thm]{Corollary}
\newtheorem{ex}{Example}
\newtheorem{defi}{Definition}
\newcommand{\pil}{\rightarrow}
\newcommand{\bpil}{\leftarrow}
\newcommand{\dpil}{\leftrightarrow}
\newcommand{\nrd}{\text{NRD}}
\newcommand{\cay}{\text{Cay}}

\begin{document}
\maketitle

\begin{abstract}
A normally regular digraph with parameters $(v,k,\lambda,\mu)$ is a
directed graph on $v$ vertices whose adjacency matrix
$A$ satisfies the equation $AA^\text{t}=k I+\lambda (A+A^\text{t})+\mu
(J-I-A-A^\text{t})$. This means that every vertex has out-degree $k$, a pair
of non-adjacent vertices have $\mu$ common out-neighbours, a pair of
vertices connected by an edge in one direction have $\lambda$ common
out-neighbours and a pair of vertices connected by edges in both
directions have $2\lambda-\mu$ common out-neighbours. We often assume
that two vertices can not be connected in both directions. 

We prove that the adjacency matrix of a normally regular digraph is
normal. 
A connected $k$-regular digraph with normal adjacency matrix is a normally regular digraph if and
only if all eigenvalues other than $k$ are on one circle in the
complex plane.
We prove a Bruck-Ryser type condition for existence and give a
combinatorial proof for a restriction excluding existence in some
cases with small values of $\lambda$.
There is a structural characterization of normally regular digraphs
with $\mu=0$ or $\mu=k$. For other values of $\mu$ we give several
constructions of normally regular digraphs. In many cases these graphs
are Cayley graphs of abelian groups and the construction is then based
on a generalization of difference sets. 
In particular, if $4t+1$, $4s+3$ and $q$ are prime powers and $r$ is not divisible by 3
we get normally regular Cayley digraphs with the following parameters
$$((4t+1)(4s+3),(4t+2)(2s+1),4st+3s+t+1),$$
$$((4s+3)(2s+1),4s+1,s,1),$$
$$(\frac{q^{2r}+q^r+1}{q^2+q+1},q^r-q,q^2,q^2+q+1)$$
and, if $q\equiv 1 \mod 3$ 
$$(\frac{1}{3}(q^2+q+1),q-1,1,3).$$

We also show connections to other combinatorial objects:
strongly regular graphs, symmetric 2-designs and association schemes.

\end{abstract}

\section{Introduction}

In this section we introduce normally regular digraphs and other basic
concepts. In Section 2 we prove that the adjacency matrices of
normally regular digraphs are normal and we give a Bruck-Ryser type
condition for existence. In Section 3 we show that complements of normally
regular digraphs are normally regular and we prove bounds on the parameters. 
In Section 4 we characterize normally
regular digraphs with $\mu=0$ or $\mu=k$. We consider eigenvalues of
normally regular digraphs in Section 5 
and show that a regular digraph with normal adjacency matrix is a normally regular digraph if and only if the non-trivial eigenvalues are on a circle in the complex plane.
In Section 6 we consider
relations to association schemes. The subject of Section 7 is
partitions of the vertex set and in particular 
group divisible digraphs, i.e.,
orientations of complete multipartite graphs.
In Section 8 we exclude existence for
some parameter sets with small $\lambda$. Section 9 describes
applications of normally regular digraphs to partitions of designs in smaller designs. 
In Section 10 we give
several constructions of normally regular digraphs, primarily
constructions as Cayley graphs.

\bigskip

The adjacency matrix of a digraph  with vertex set
$\{x_1,\ldots x_v\}$ is a $v\times v$ matrix $A$ in which the
$(i,j)$-entry is
$$A_{ij}=
\begin{cases}
1 & \text{if there is an edge directed from $x_i$ to $x_j$ }\\
0 & \text{otherwise.}
\end{cases}$$
Thus any square $\{0,1\}$-matrix is the adjacency matrix of a digraph if
and only if all its diagonal entries are 0.
In this paper we consider such matrices that satisfy an equation
involving $AA^\text{t}$.
The $(i,j)$ entry of $AA^\text{t}$ (respectively $A^\text{t}A$)
is the number of common out-neighbours (respectively in-neighbours) of
$x_i$ and $x_j$.

We say that a digraph is normal if its adjacency matrix $A$ is normal,
i.e., if $AA^\text{t}=A^\text{t}A$. It follows that a 
digraph is normal if and only if for any two (not necessarily
distinct) vertices $x$ and $y$ the number of common out-neighbours of
$x$ and $y$ is equal to the number of common in-neighbours of $x$ and
$y$.

We will use the notation $x\pil y$ if there is an edge directed from
$x$ to $y$ (and possibly also an edge from $y$ to $x$). If $x\pil y$ then
we say that $x$ dominates $y$.
We write $x\dpil y$
if $x\pil y$ and $x\bpil y$, and identify these two directed edges
with an undirected edge.

The set $\{y\mid x\pil y\}$ of out-neighbours of a vertex $x$ is
denoted by $x^+$. Similarly $x^-$ denotes the set of in-neighbours. 
$d^+(x)=|x^+|$ and $d^-(x)=|x^-|$ denotes the out-degree and in-degree
of $x$, respectively. 

We will now introduce normally regular digraphs. We first give a
matrix free definition.

\begin{defi} \label{nrd-def}
A normally regular digraph with parameters $(v,k,\lambda,\mu)$,
also denoted by $\nrd(v,k,\lambda,\mu)$,
is a directed graph on $v$ vertices so that
\begin{itemize}
\item every vertex has out-degree $k$
\item any pair of non-adjacent vertices have exactly $\mu$ common
  out-neighbours,
\item any pair of vertices $x,y$ such that exactly one of the edges $x\pil
  y$ or $x\bpil y$ is present have exactly $\lambda$ common out-neighbours,
\item any pair of vertices $x,y$ such that $x\dpil y$ have exactly
  $2\lambda-\mu$ common out-neighbours.
\end{itemize}
A normally regular digraph is said to be asymmetric if there is no
pair $x,y$ so that $x\dpil y$.
\end{defi}

\newpage
This definition may be stated in terms of the adjacency matrix.

\begin{prop} \label{prop:matdef}
A $v\times v$ \ $\{0,1\}$-matrix $A$ is the adjacency matrix of a normally
regular digraph if and only if every diagonal entry is 0 and
$$AA^\text{t}=kI+\lambda(A+A^\text{t})+\mu (J-I-A-A^\text{t}),$$
where $I$ is the identity matrix and $J$ is the matrix in which all
entries are 1.

This normally regular digraph is asymmetric if and only if
$A+A^\text{t}$ is a $\{0,1\}$ matrix.
\end{prop}

The author first intended  to study only
asymmetric normally regular digraphs. 
However, most of the results  
hold in the general case, so we will usually not assume that
graphs are asymmetric, but for connections to association schemes and
similar results we need to assume that the graph is asymmetric. 

Asymmetric normally regular digraphs with $\mu=\lambda$ have been
studied by Ito~\cite{Ito88}, \cite{Ito89a}, \cite{Ito89b},
\cite{Ito89c} and \cite{Ito90}, and also by Ionin and
Kharaghani~\cite{IKharaghani}.

Fossorier, Je\v{z}ek, Nation and Pogel~\cite{ordinary} introduced what
they call ordinary graphs. Their definition is similar to our
Definition~\ref{nrd-def}, but the number of common out-neighbours (and
common in-neighbours) of $x$ and $y$ in the three cases is $a$, $b$
and $c$, respectively. They do not assume that $c=2b-a$ (although this
is satisfied in some of their results). Note that the equation
$c=2b-a$ is essential for the formulation of the definition of a normally
regular digraph as a matrix equation, and thus it is
essential for the theory.

U. Ott~\cite{Ott} considered Cayley graph construction from
``generalized difference set'' that leads to normally regular digraphs.

Another variation of normally regular digraphs is Deza digraphs.
A regular digraph is said to be a Deza digraph if the number common
out-neighbours of two vertices is either $b$ or $c$, for some
constants $b$ and $c$, but it need not depend on whether the vertices
are adjacent or not. Deza digraphs have been studied by Wang and
Feng~\cite{WangFeng}.  

Many constructions of normally regular digraphs uses Cayley graphs of a
group. 
Let $G$ be a group and let $S$ be a subset of $G$ not containing the
group identity. Then the Cayley
graph $\cay(G,S)$ is the graph whose vertices are the elements of $G$
and with edge set
$$\{x\pil y \mid x^{-1}y\in S\}.$$ 
Let $S^{(-1)}=\{s^{-1}\mid s\in S\}$. Then $\cay(G,S)$ is undirected
if $S^{(-1)}=S$ and asymmetric if  $S^{(-1)}\cap S=\emptyset$.
$\cay(G,S)$ is a normally regular digraph if for every $g\in G, g\neq
1$, the number of pairs $(x,y)\in S\times S$ satisfying $yx^{-1}=g$
is $\mu$ if $g\notin
S\cup S^{(-1)}$, $\lambda$ if $g$ is in exactly one the sets $S,
S^{(-1)}$ and $2\lambda-\mu$ if $g\in S\cap S^{(-1)}$. 

In~\cite{nrd2} we prove the following multiplier theorem.

\begin{thm} \label{multthm}
Suppose that $G$ is an abelian group and that $\cay(G,S)$ is an
$\nrd(v,k,\lambda,\mu)$.
Let $w$ be the smallest positive
number so that for every $g\in G$ the order of $g$ divides $w$.
Let $m$ be a natural number relatively prime to $v$, so that $m$ divides
$\eta=k-\mu+(\lambda-\mu)^2$ and let $t$ be relatively prime to $v$.

Suppose that for every prime $p$ dividing $m$ there exist an integer
$f$ so that $t\equiv p^f \mod w$.
If either $m>\mu\geq\lambda+2$ or $m>2\lambda-\mu$, $\lambda>\mu$ then
(in additive notation) $tS=S$.
\end{thm}

Furthermore, in~\cite{nrd3} we enumerate small normally regular
digraphs and prove some results related these graphs.
In \cite{JJKS}, group divisible normally regular digraphs, i.e, the digraphs considered in Section~\ref{mu=0} and Section~\ref{groupdivisible}, are investigated.

It is well-known (see Godsil and Royle~\cite{GodRoy}) that a strongly regular graph with parameters
$(v,k,\lambda,\mu)$ is an undirected graph with $v$ vertices in which
\begin{itemize}
\item every vertex has degree $k$
\item any pair of adjacent vertices have exactly $\lambda$ common
  neighbours
\item any pair of non-adjacent vertices have exactly $\mu$ common
  neighbours.
\end{itemize}
Equivalently, a strongly regular graph is a graph whose adjacency
matrix $A$ satisfies
$$A^2=kI+\lambda A +\mu (J-I-A)\quad \text{and} \quad AJ=JA=kJ.$$

Thus any normally regular digraph where all edges are undirected
(i.e., $x\pil y$ if and only if $x\bpil y$) is a strongly regular
graph. Note however that we use $\lambda$ in a different meaning. For
a normally regular digraph we will use $\lambda_2=2\lambda-\mu$ to
denote the number of common out-neighbours of a pair of vertices joined by
two edges. 

In the theory of normally regular digraphs we will require that
$\lambda$ is an integer and thus $\mu$ and $\lambda_2$ are congruent
modulo 2. Thus not every strongly regular graph is a normally regular
digraph. 

A strongly regular with $(v,k,\lambda,\mu)=(4\mu+1,2\mu,\mu-1,\mu)$ is
called a conference graph. Since $\lambda$ and $\mu$ have different
parity a conference is not a normally regular digraph, but it will used
in some constructions. The most important construction of conference
graphs are the Paley graphs which are constructed as follows. Let $F$
be a field of $q$ elements, $q\equiv 1 \mod 4$ and let $Q$ be the
non-zero squares in $F$. Then the Cayley graph of the additive group
$\cay(F,Q)$ is a conference graph, see~\cite{GodRoy}.

There are some directed analogues of strongly regular graphs other
than normally regular digraphs. Duval~\cite{Duval} introduced directed
strongly regular graphs which have adjacency matrix $A$ satisfying
$$A^2=tI+\lambda A+\mu (J-I-A)\text{\ \ and \ \ } AJ=JA=kJ.$$
Many proof techniques from strongly regular graphs, especially the use
of eigenvalues, are more easily applied to directed strongly regular
graphs than to normally regular graphs, see~\cite{Duval} or~\cite{korea}.

Another well-known combinatorial structure to which normally regular
digraphs are related are 2-designs (or Balanced Incomplete Block Designs). 
A $2-(v,k,\lambda)$ design is an incidence structure with $\{0,1\}$
incidence matrix $N$ of size $v\times b$, $b=\frac{\lambda
  v(v-1)}{k(k-1)}$ satisfying 
\begin{equation} \label{design-def}
NN^\text{t}=(k-\lambda)I +\lambda J. 
\end{equation}
A 2-design is said to be symmetric if $b=v$. The parameter $k-\lambda$ is
called the order of the symmetric design.
For information on design theory, see Beth, Jungnickel and Lenz~\cite{Beth}.

Let $A$ be the adjacency matrix of a normally regular digraph.
If $\mu=\lambda$ then $A$ is incidence matrix of symmetric 2-design.
If $\mu=\lambda+1$ then $A+I$ is incidence matrix of a symmetric 2-design.
In this paper we will often assume that $\mu\notin\{\lambda,\lambda+1\}$.

A tournament is a digraph with the property that for any two distinct
vertices $x$ and $y$ exactly one of the edges $x\pil y$ or $y\pil x$
is present.

We will need the following property of regular tournaments in
Section~\ref{mu=0}.

\begin{lemma} [Rowlinson~\cite{Rowlinson}] \label{normaltournament}
A tournament is normal if and only if it is regular.
\end{lemma}
\proof 
If $A$ is an adjacency matrix of a regular tournament, i.e.,
$AJ=JA=kJ$, for some number $k$ then, since $A^\text{t}=J-I-A$,
$AA^\text{t}= AJ-A-A^2=JA-A-A^2=A^\text{t}A$. 

Conversely, in a normal digraph every vertex has the same in-degree and
out-degree. \qed

\medskip

If a tournament is a normally regular digraph then it is called a doubly
regular tournament. It satisfies $k=2\lambda+1$.
$\mu$ is arbitrary, so we may take $\mu=\lambda$. Such tournaments
are also called homogenous tournaments by Kotzig~\cite{Kotzig},
and Ito~\cite{Ito90} used the term Hadamard tournaments, as these tournaments
are equivalent to skew Hadamard matrices of order $v+1$ (see Reid and
Brown~\cite{ReidBrown}). 

Thus it is possible that doubly regular tournaments of order $v$ exists
for all $v\equiv 3 \mod 4$.

The most important construction of a doubly regular tournament is the
Paley tournament which is constructed as follows. Let $F$ be a field
of $q$ elements $q\equiv 3 \mod 4$ and let $Q$ be the non-zero
squares in $F$. Then the Cayley graph $\cay(F,Q)$ is a doubly regular
tournament.

We conclude this section with two small asymmetric normally regular
digraphs.

\begin{ex} \label{ex-intro}
Let $Q=\{1,i,j,k,-1,-i,-j,-k\}$ be the quaternion group. Then
$\cay(Q,\{i,j,k\})$ is an NRD(8,3,1,0) with the following adjacency matrix
$$\begin{pmatrix}
0&1&1&1&0&0&0&0\\
0&0&0&1&1&0&1&0\\
0&1&0&0&1&0&0&1\\
0&0&1&0&1&1&0&0\\
0&0&0&0&0&1&1&1\\
1&0&1&0&0&0&0&1\\
1&0&0&1&0&1&0&0\\
1&1&0&0&0&0&1&0
\end{pmatrix}$$
This
is the smallest non-trivial normally regular digraph with $\mu=0$. 
Normally regular digraphs with $\mu=0$ or $\mu=k$ are characterized in
Section~\ref{mu=0}. 
\end{ex}

\begin{ex}
$\cay(\mathbb{Z}_{19},\{1,4,6,7,9,11\})$ is an NRD(19,6,1,3). This is the
smallest normally regular digraph with
$\mu\notin\{k,0,\lambda,\lambda+1\}$. It belongs to an infinite
family constructed in Theorem~\ref{k13}. This digraph is asymmetric
and in fact $\lambda_2=2\lambda-\mu$ is negative.
\end{ex}

\section{Matrix equations}

It is convenient to introduce two further parameters of a normally
regular digraph:
$$\eta=k-\mu+(\mu-\lambda)^2$$
and
$$\rho=k+\mu-\lambda.$$
The parameter $\eta$ will play a role similar to that of the order of
a symmetric design. Ma~\cite{Ma} uses the parameter $\Delta=4\eta$ in
the study of strongly regular graphs. The factor 4 is necessary in order
get an integer for a general strongly regular graph.

The matrix equation in Proposition~\ref{prop:matdef} is equivalent the
to following equation.
\begin{equation} \label{design-eq}
(A+(\mu-\lambda)I)(A+(\mu-\lambda)I)^\text{t}=\eta I +\mu J.
\end{equation}
Thus for $B=(A+(\mu-\lambda)I)$ we have
$$BB^\text{t}=\eta I +\mu J,$$
and since $AJ=kJ$ (every vertex has out-degree $k$),
$$BJ=\rho J.$$

We will now prove that a normally regular digraph is normal. 
The following lemma is a generalization of a proof of the
fact that the dual of a symmetric 2-design is also a 2-design,
see~\cite{Beth}.  

\begin{lemma}
Suppose that $B$ is a non-singular $v\times v$ matrix so that
$BB^\text{t}=\eta I +\mu J$ and $BJ=\rho J$ for some constants $\rho,
\eta, \mu$. Then $B$ is normal and $\mu v=\rho^2-\eta$.
\end{lemma}

\proof
From $BJ=\rho J$ we get $\rho^{-1} J= B^{-1}J$ and
\begin{equation}  \label{eq:B}
B^\text{t}=B^{-1}(BB^\text{t})=B^{-1}(\eta I+\mu J)=\eta
B^{-1}+\mu\rho^{-1} J
\end{equation}
Using that $J$ is symmetric, we get from this
$$\rho J= (BJ)^\text{t}=JB^\text{t}=\eta JB^{-1}+\mu\rho^{-1} J^2=\eta
JB^{-1}+\mu\rho^{-1}v J.$$
This implies that
$$JB^{-1}=\frac{\rho-\mu \rho^{-1}v}{\eta} J,$$
and so
$$vJ=J^2=(JB^{-1})(BJ)=\frac{\rho-\mu\rho^{-1}v}{\eta}\rho v J.$$
Thus
\begin{equation}
\frac{\rho-\mu\rho^{-1}v}{\eta}=\rho^{-1},   \label{eq:veq1}
\end{equation}
and $JB^{-1}=\rho^{-1}J$ or $\rho J=JB$.
Now equation~\ref{eq:B} implies
$$B^\text{t}B=\eta I+\mu\rho^{-1} JB=\eta I+\mu J=BB^\text{t}.$$
Rewriting equation~\ref{eq:veq1} we get
$\mu v=\rho^2-\eta$. \qed

\begin{cor} \label{NRDisnormal}
Every normally regular digraph is normal.
\end{cor}

\proof
Let $A$ be the adjacency of a normally regular digraph and let $B=A+(\mu-\lambda)I$.
Then $BB^\text{t}=\eta I+\mu J$. Suppose first
that $B$ is singular. Then one of the eigenvalues of $\eta I+\mu J$ is
zero: $\eta=0$ or $\eta+\mu v=0$. Since $\mu,v\geq 0$ this is possible
only when $\eta=k-\mu+(\mu-\lambda)^2$ is $0$.
As $k+(\mu-\lambda)^2\geq k \geq \mu$,
$\mu=k+(\mu-\lambda)^2$ implies $k=\mu=\lambda$.
This implies that $k=0$. Since a graph with
no edges is normal, we may thus assume that $B$ is non-singular, and
the result follows from the lemma.\qed

\medskip
It follows that a normally regular digraph is both normal and regular,
i.e., every vertex has in-degree $k$ and out-degree $k$. And the
number of common in-neighbours of distinct vertices $x$ and $y$ is
$$\begin{cases}
\mu & \text{\ if $x$ and $y$ are non-adjacent,}\\
\lambda & \text{\ if either\ } x\pil y \text{\ or\ } y\pil x, \text{\
  but not both,}\\
2\lambda-\mu & \text{\ if $x\dpil y$}.
\end{cases}$$

\begin{cor}   \label{param-eq}
The parameters of a normally regular digraph satisfy
\begin{equation}
\mu v=\rho^2-\eta.  \label{eq:veq2}
\end{equation}
\end{cor}

This equation is equivalent to the following
\begin{equation}
2k\lambda+(v-2k-1)\mu=k^2-k  \label{param-eq1}
\end{equation}

This equation may also be obtained by counting in two ways the number of
triples $(x,y,z)$ of vertices so that $x\pil y\bpil z$ using the
definition of a normally regular digraph and the fact that every
vertex has in-degree $k$. 

From the theory of symmetric 2-designs we also have the Bruck-Ryser
type condition. It is based on the following general lemma from Beth,
Jungnickel and Lenz~\cite{Beth}

\begin{lemma} \label{Bethlemma}
Suppose that $N$ is a rational $v\times v$ matrix satisfying the
equation
$$NN^\text{t}=(a-b)I+bJ$$
for some integers $a>b$ and $v$ odd. Then the equation
$$x^2=(a-b)y^2+(-1)^{(v-1)/2}bz^2$$
has a solution $(x,y,z)\in \mathbb{Z}^3\setminus\{(0,0,0)\}$.
\end{lemma}

For normally regular digraphs we have the following.

\begin{thm} \label{BRC}
Suppose that there exist an $\nrd(v,k,\lambda,\mu)$.
\begin{itemize}
\item If $v$ is even then $\eta=k-\mu+(\mu-\lambda)^2$ is a square.
\item If $v\equiv 1$ {\em (mod 4)} then the Diophantine equation $x^2-\mu
  y^2=\eta z^2$ has an integer solution such that
  $x$, $y$, and $z$ are not all zero.
\item If $v\equiv 3$ {\em (mod 4)} then the Diophantine equation $x^2+\mu
  y^2=\eta z^2$ has an integer solution such that
  $x$, $y$, and $z$ are not all zero.
\end{itemize}
\end{thm}

\medskip
\noindent {\bf Proof\ } It follows from equation~\ref{design-eq} that
the determinant of $\eta I+\mu J$ is a square. The
eigenvalues of this matrix are
$\eta+\mu v=\rho^2$ with multiplicity 1 and
$\eta$ with multiplicity $v-1$. For the equality we
used equation~\ref{eq:veq2}. Thus the result follows when $v$ is even.

For $v$ odd, the theorem follows from
equation~\ref{design-eq} and the above lemma. \qed

\section{Complementary graphs and the parameters}

The complement of a graph with adjacency matrix $A$ is the graph with
adjacency matrix $J-I-A$. The following theorem is proved by an easy
computation.

\begin{thm} \label{complement}
Let $A$ the adjacency matrix of a normally regular
digraph with parameters $(v,k,\lambda,\mu)$. Then $J-I-A$ is the
adjacency matrix of a normally regular digraph with
parameters 
$$(\overline{v},\overline{k},\overline{\lambda},\overline{\mu})=
(v,v-k-1,v-2k+\lambda-1,v-2k+2\lambda-\mu).$$
\end{thm}

Note that $\overline{\eta}=\overline{k}-\overline{\mu}+
(\overline{\mu}-\overline{\lambda})^2=\eta$.

Two important cases, $\mu=k$ and $\mu=0$, are considered in the next
section. They are complementary.
\begin{cor} \label{my=k0}
A normally regular digraph satisfies $\mu=k$ if and only
if the complementary normally regular digraph satisfies
$\mu=0$.
\end{cor}

\proof 
If $\mu=0$ then it follows from equation~\ref{param-eq1} that
$2\lambda=k-1$. And then $\overline{\mu}=v-2k+2\lambda-\mu=v-k-1=
\overline{k}$.

If $\mu=k$  then it follows from equation~\ref{param-eq1} that
$v=3k-2\lambda$. And then $\overline{\mu}=v-2k+2\lambda-\mu=0$. \qed

\medskip
We will now consider upper and lower bounds on the parameters $\mu$,
$\lambda$ and $\lambda_2$.
There exists normally regular digraphs for which
$\lambda_2=2\lambda-\mu<0$. But in that case there can not be any
undirected edges and so the digraph is asymmetric. 
Note that Theorem~\ref{complement} is still valid in this case.

\begin{lemma} \label{lemma:paramulighed}
The paramaters of an asymmetric normally regular digraph with $k\geq
1$ satisfy the following restriction: 
$$k\geq 2\lambda+1.$$
\end{lemma}

\proof
The number of edges in the subgraph spanned by the set $x^+$ of
out-neighbours of a vertex $x$ is $k\lambda\leq \binom{k}{2}$. Thus
$2\lambda\leq k-1$.\qed

\medskip
If a normally regular digraph is a tournament then $\mu$ can be chosen
arbitrarily and if it is a complete undirected graph $K_v$ then $\mu$
and $\lambda$ can be chosen arbitrarily so that $\lambda_2=2\lambda-\mu=v-2$.
In all other cases we have $0\leq \mu\leq k$.

\begin{prop} \label{mu-interval}
Suppose there exists an $\nrd(v,k,\lambda,\mu)$ which is not a
tournament or a complete graph. Then
$$0\leq\mu\leq k,$$
and 
$$\lambda_2=2\lambda-\mu\leq k-1,$$
with equality if and only if $\mu=0$.
\end{prop}

\proof
Suppose that $\mu<0$. Then the digraph does not have any
pair of non-adjacent vertices. Since it is not a complete graph,
$v-2k-1>-k$. As the digraph is not a tournament, there
exist undirected edges and $2\lambda-\mu=\lambda_2\leq k-1$.
From equation~\ref{param-eq1} we have 
$$k^2-k=2k\lambda+(v-2k-1)\mu< 2k\lambda-k\mu\leq k(k-1),$$
a contradiction. Thus $\mu\geq 0$.

If $\mu=0$ then by equation~\ref{param-eq1}, $2\lambda=k-1$ and so
$\lambda_2=k-1$. 
If the digraph has an undirected edge then clearly $\lambda_2\leq
k-1$. If $x\dpil y$ and $x$ and $y$ have $\lambda_2=k-1$ common
out-neighbours then in the complementary graph they have
$v-k-1=\overline{k}$ common out-neighbours. Thus
$\overline{\mu}=\overline{k}$ and by Corollary~\ref{my=k0}, $\mu=0$.
So suppose that the digraph is asymmetric. Then by
Lemma~\ref{lemma:paramulighed}, $2\lambda\leq k-1$ and so
$2\lambda-\mu\leq k-1$, with equality only if $\mu=0$. 

Suppose now that $\mu>k$. Let
$(v,\overline{k},\overline{\lambda},\overline{\mu})$ be the parameters
of the complementary normally regular digraph. Then by
Theorem~\ref{complement} 
$$v-2\overline{k}+2\overline{\lambda}-\overline{\mu}=\mu >k=v-\overline{k}-1,$$
and so $2\overline{\lambda}-\overline{\mu}>\overline{k}-1$, a
contradiction. \qed

\section{$\mu=0$ or $\mu=k$\label{mu=0}}

By Theorem~\ref{my=k0}, normally regular digraphs with $\mu=0$ and
$\mu=k$ are complements of each other. We will therefore characterize
normally regular digraphs with $\mu=0$ and then get the case $\mu=k$
as a corollary.

\subsection{$\mu=0$}

We will first characterize {\em asymmetric} normally regular digraphs with
$\mu=0$ and then generalize to digraphs with undirected edges.

A normally regular digraph with $\mu=0$ need not be connected. However, each connected
component will be a normally regular digraph with the same value of $k$ and $\lambda$.
Thus we will only consider normally regular digraphs whose underlying undirected graph is
connected. As each vertex has equal in- and out-degree this implies
that the digraph is strongly connected. Thus there is a directed path
from any vertex to any other vertex. A normally regular digraph with $\mu=0$ may be a
doubly regular tournament. Another possibility is that
$k=1$ and the digraph is a directed cycle.

Let $T$ be a tournament with adjacency $A$. Then $\mathcal{D}(T)$
denotes the digraph with adjacency matrix
$$\left (
\begin{array}{cccc}
0 & 1 \ \ldots \  1 & 0 & 0 \ \ldots \  0\\
0 &            & 1 &         \\
\vdots & A     &\vdots & A^\text{t}  \\
0 &  & 1 &  \\
0 & 0 \ \ldots \  0 & 0 & 1 \ \ldots \  1\\
1 &            & 0 &         \\
\vdots & A^\text{t}     &\vdots & A  \\
1 &  & 0 &  \\
\end{array}\right ).
$$
Thus if $T$ is a tournament with only one vertex then $\mathcal{D}(T)$
is a directed cycle of length 4. In this section we consider a
tournament with one vertex to be doubly regular.

\begin{thm} \label{asym-my=0}
A connected digraph is an asymmetric normally regular digraph with $\mu=0$ if and only if either
\begin{enumerate}
\item it is a directed cycle of length at least 5
\item it is a doubly regular tournament or
\item it is isomorphic to $\mathcal{D}(T)$ for some doubly regular
  tournament $T$.
\end{enumerate}
\end{thm}

\proof
Suppose that $G$ is a connected asymmetric normally regular digraph with $\mu=0$, $k\geq
2$ and that $G$ is not a tournament.

As $\mu=0$ we get from 
equation~\ref{param-eq1} that
$\lambda=\frac{k-1}{2}$. Let $x$ be a vertex
of $G$. Then every vertex in $x^+$ has out-degree $\lambda$ in this
subgraph and thus in-degree at most $k-1-\lambda=\lambda$. It follows
that $x^+$ is a regular tournament. Similarly, $x^-$ is a regular tournament.

Since $G$ is strongly connected and it is not a tournament, there
exist a vertex $y\in G$ so that $x$ and $y$ are non-adjacent and there
is a path from $x$ to $y$ in $G$. We may choose $y$ so that the
(directed) distance from $x$ to $y$ is minimal, i.~e. $y$ is dominated
by a vertex in $x^+$ or in $x^-$. Since $x$ and $y$ are non-adjacent
and $\mu=0$, $y$ is not dominated by any vertex in $x^-$, and
similarly $y$ does not dominate any vertex in $x^+$. Thus $y$ is
dominated by a vertex, say $v$, in $x^+$.
Suppose there is a vertex $w$ in $x^+$
that does not dominate $y$. Since $x^+$ is a regular tournament it is
strongly connected, so there is a directed path from $v$ to $w$ in
$x^+$. On this path there are vertices $u$ and $u'$ so that
$u\rightarrow u'$, $u\rightarrow y$ but $u'$ does not dominate $y$.
This is a contradiction to $\mu=0$. Thus every vertex in $x^+$
dominates $y$. If another vertex $y'$, non-adjacent to $x$ was
dominated by a vertex in $x^+$, it would be dominated by every vertex
in $x^+$ and so $y$ and $y'$ have $k$ common in-neighbours, a
contradiction. Thus every vertex in 
$x^+$ dominates $\lambda $ vertices in $x^-$.

Now a vertex in $x^-$ dominated by a vertex in $x^+$ (which dominates
$y$) must be adjacent to $y$, as $\mu=0$. As above, $y$ then dominates
every vertex of $x^-$, and every vertex in $x^-$ is dominated by
$\lambda$ vertices in $x^+$. Also every vertex in $x^-$ dominates
exactly $\lambda$ vertices in $x^+$. Thus $V(G)=\{x,y\}\cup x^+ \cup
x^-$. Furthermore there is an enumeration of vertices
$x^+=\{v_1,\ldots,v_n\}$ and $x^-=\{v_1',\ldots,v_n'\}$ such that
$v_i'$ is the unique vertex non-adjacent to $v_i$ and vice versa.

If $v_i\rightarrow v_j$ then, since no vertex dominates both $v_j$ and
$v_j'$, $v_j'\rightarrow v_i$. Similarly $v_j\rightarrow v_i'$ and
$v_i'\rightarrow v_j'$. Thus the mapping $v_i\mapsto v_i'$ is an isomorphism.

We also see that $v_\ell'$ is a common out-neighbour of $v_i$ and $v_j$  if
and only if $v_\ell$ is a common in-neighbour of $v_i$ and $v_j$.
Thus the number of vertices in $x^+$
dominating $v_i$ and $v_j$ plus the number of vertices in $x^+$
dominated by $v_i$ and $v_j$ is $\lambda-1$. But since $x^+$ is
a regular tournament it is normal (by Lemma~\ref{normaltournament})
and thus these two numbers are both
equal to $\frac{\lambda-1}{2}$ and so $\lambda$ is odd, and $x^+$ is
a doubly regular tournament, $\nrd(k,\lambda,\frac{\lambda-1}{2},\cdot)$.

If on the other hand $G$ is a doubly regular tournament with degree
$k=2\lambda+1$ and with vertex-set $\{x_1,\ldots,x_n\}$, $n=2k+1$,
then we may construct a graph with vertex-set
$\{v_0,\ldots,v_n,v_0',\ldots, v_n'\}$ and edges
$$v_0\rightarrow v_i\rightarrow v_0'\rightarrow v_i'\rightarrow v_0,
\mbox{\ for \ } 1\leq i\leq n$$
and
$$v_i\rightarrow v_j\rightarrow v_i'\rightarrow v_j'\rightarrow v_i
\mbox{\ if\ } x_i\rightarrow x_j \mbox{\ in\ } G, \mbox{\ for \ }
1\leq i,j\leq n.$$
It is easy to verify that this new graph is an $\nrd(2n+2,n,k,0)$.
\qed

\medskip
The smallest non-trivial example of the type of normally regular digraphs
mentioned as possibility 3 is a Cayley graph of the
quaternion group of order 8 (see
Example~\ref{ex-intro}). In~\cite{switch} and~\cite{JJKS}, it is 
investigated when normally regular digraphs of this type are Cayley
graphs or vertex transitive.

\bigskip

We will now characterize normally regular digraphs with $\mu=0$ and with
undirected edges. We need a definition to describe the digraphs.
Let $G$ be a digraph with vertices $x_1,\ldots,x_n$. Then we denote by
$\mathcal{K}_s(G)$ the digraph with vertex set partitioned in sets
$V_1,\ldots,V_n$ of size $s$ where each $V_i$ induce an complete
undirected graph and furthermore for $y\in V_i$ and $z\in V_j$, $y\pil
z$ if and only if $x_i\pil x_j$ in $G$. If $B$ is an adjacency matrix
of $G$ then an adjacency matrix of $\mathcal{K}_s(G)$ can be expressed
using Kronecker products of matrices (see Hall~\cite{Hall}) as follows
$B\otimes J_s+I_n\otimes(J_s-I_s)=(B+I)\otimes J_s-I_{ns}$.

\begin{thm} \label{generel-my=0}
Let $\Gamma$ be a connected  normally regular
digraph with parameters $(v,k,\lambda,0)$, i.e., $\mu=0$.
Then for some number $s$ there is an asymmetric normally regular
digraph $\Gamma'$ with 
parameters $(\frac{v}{s},\frac{k-s+1}{s},\frac{\lambda-s+1}{s},0)$ 
so that $\Gamma$ is isomorphic to $\mathcal{K}_s(\Gamma')$.

Conversely, if $\Gamma'$ is an asymmetric normally regular digraph with
parameters $(v,k,\lambda,0)$ then $\mathcal{K}_s(\Gamma')$ is a normally
regular digraph with parameters $(sv,sk+s-1,s\lambda+s-1,0)$.
\end{thm}

\proof
Consider a connected normally regular digraph with $\mu=0$. We have
that $k=2\lambda+1$. 
Then the number of common out-neighbours of $x$ and $y$, where
$x\leftrightarrow y$, is $2\lambda-\mu=k-1$.

Thus if $x\leftrightarrow y$ then $x$ and $y$ have exactly the same
set of out-neighbours (and the same set of in-neighbours) other than
$y$ and $x$. In particular, if $x\leftrightarrow y\leftrightarrow z$
then $x\leftrightarrow z$.

It follows that the vertex set is partitioned in sets
$V_1,\ldots,V_m$, so that each $V_i$ spans a complete subgraph and
there are no undirected edges joining $V_i$ and $V_j$ for $i\neq
j$. If $x\pil y$ for some $x\in V_i$ and $y\in V_j$ then $x\pil y$ for
every $x\in V_i$ and $y\in V_j$.

Choose $i$ so that $|V_i|\geq |V_j|$, for all $j$ and let $s=|V_i|$.
Let $V_i^+$ denote the set of out-neighbours outside $V_i$ of vertices
in $V_i$.
Then $V_i^+=V_{i_1}\cup\ldots\cup V_{i_\ell}$ for some $i_1,\ldots,i_\ell$
and $V_i^+$ has size
$k-(s-1)=2\lambda+2-s$. In the subgraph spanned by
$V_i^+$ every vertex has out-degree $\lambda$. The average
in-degree is also $\lambda$. Thus the average number of undirected edges
incident with a vertex is at least $2\lambda-(2\lambda+2-s-1)=s-1$. By the
maximality of $|V_i|=s$, no vertex is incident with more than $s-1$
undirected edges and so $|V_{i_1}|=\ldots=|V_{i_\ell}|=s$. Since the graph
is connected, repeated use of this argument shows that
$|V_1|=\ldots=|V_m|$.
Consider a graph $\Gamma'$ with vertices $x_1,\ldots,x_m$ and edges
$x_i\pil x_j$ if $V_i\pil V_j$.Then $\Gamma'$ is a normally regular
digraph with parameters
$(\frac{v}{s},\frac{k-s+1}{s},\frac{\lambda-s+1}{s},0)$, and $\Gamma$
is isomorphic to $\mathcal{K}_s(\Gamma')$.

\qed

\newpage
\begin{ex} \label{Q16}
If $\Gamma'$ in this proof is a normally regular digraph with
parameters $(8t+8,4t+3,2t+1,0)$ then the parameters of $\Gamma$ are
$((8t+8)s,(4t+4)s-1,(2t+2)s-1,0)$. Thus a normally
regular digraph with parameters $(16,7,3,0)$ may appear with
$(s,t)=(2,0)$ or with $(s,t)=(1,1)$ and so $s$ can not be determined
from the parameters.
\end{ex}

\subsection{$\mu=k$}
\begin{thm} \label{myk}
A digraph $G$ is an asymmetric normally regular digraph with $\mu=k$ if and only if there
is a number 
$s$ so that $G$ is obtained from a doubly-regular tournament by
replacing each vertex $x$ by a set $V_x$ of $s$ new vertices such
that if $x\rightarrow y$ in the tournament then $u\rightarrow w$ for
every $u\in V_x$ and $w\in V_y$. Then $s=k-2\lambda=v-2k$.
\end{thm}

In other words a graph is an asymmetric normally regular digraph with $\mu=k$ if and only if it has an
adjacency matrix which is the Kronecker product of an
adjacency matrix of a doubly regular tournament and $J_s$

\medskip
\proof
If $G$ is an asymmetric normally regular digraph with $\mu=k$ then the
complement $\overline{G}$ of $G$ is a connected normally regular digraph with
$\mu=0$ and with no pair of non-adjacent vertices. Then $\overline{G}$
is constructed as in Theorem~\ref{generel-my=0} from an asymmetric
normally regular digraph with $\mu=0$ and with no pair of non-adjacent
vertices. By Theorem~\ref{asym-my=0} this is a doubly regular
tournament. \qed

\section{Eigenvalues}

If $A$ is the adjacency matrix of a normally regular digraph then we have the spectral decomposition for normal matrices
$$A=\sum_\theta \theta E_\theta,$$
where the sum is over the eigenvalues of $A$, and $E_\theta$ is the matrix of the orthogonal projection on the corresponding eigenspace. As $A$ is real, the adjoint matrix (i.e., the complex conjugate of the transposed matrix) is $A^\text{t}$ and since orthogonal projections are self-adjoint we get
$$A^\text{t}=\sum_\theta \overline{\theta} E_\theta.$$
In particular, $x\in \mathbb{C}^v$ is an eigenvector of $A$ with eigenvalue $\theta$ if and only $x$ is an eigenvector of $A^\text{t}$ with eigenvalue $\overline{\theta}$.

In general it is not possible to compute the eigenvalues of a normally
regular digraph from its parameters. We only know that the degree $k$
is an eigenvalue and (if the graph is connected then) it has
multiplicity 1. We now show that all other eigenvalues lie on a circle
in the complex plane with centre $\lambda-\mu$ and radius
$\sqrt{\eta}$.

\begin{thm} \label{thm:eigenval}
Suppose that $\theta\neq k$ is an eigenvalue of an
$\nrd(v,k,\lambda,\mu)$.

Then
\begin{equation} \label{eq:eigenvals}
|\theta-(\lambda-\mu)|=\sqrt{\eta}.
\end{equation}

\end{thm}

\proof
Let $A$ be the adjacency matrix of an $\nrd(v,k,\lambda,\mu)$.
Let
$x\in \mathbb{C}^v$ be an eigenvector for $A$ with eigenvalue
$\theta$.
Then $x$ is eigenvector for $A^\text{t}$
with eigenvalue ${\overline{\theta}}$.

Thus $(\theta+\mu-\lambda)({\overline{\theta}}+\mu-\lambda)$ is an
eigenvalue of $(A+(\mu-\lambda)I)(A+(\mu-\lambda)I)^\text{t}=\eta
I+\mu J$.

If $\theta=k$ then all entries of $x$ are equal and
$$(k+\mu-\lambda)^2=\eta+\mu v$$
(this is in fact equation~\ref{eq:veq2}).

If $\theta\neq k$ then
$$(\theta+\mu-\lambda)({\overline{\theta}}+\mu-\lambda)=\eta$$
or
$$|\theta-(\lambda-\mu)|=\sqrt{\eta}.$$
\qed

We now show that equation~\ref{eq:eigenvals} characterizes normally
regular digraphs. This theorem generalizes the well-known result that
a connected regular  undirected graph with exactly three eigenvalues
is strongly regular, see~\cite{GodRoy}.

\begin{thm}
Suppose that $G$ is a connected $k$-regular directed graph with a normal
adjacency matrix $A$. If there exist real numbers $a$ and $b$ so that
every eigenvalue $\theta\neq k$ satisfies $|\theta-a|=b$ then $G$ is a
normally regular digraph with $\lambda=a+\frac{(k-a)^2-b^2}{v}$ and
$\mu=\frac{(k-a)^2-b^2}{v}$, where $v$ is the number of vertices, or
else $G$ is a strongly regular graph.
\end{thm}

\proof
We can write $A=\sum_{i=1}^m\theta_iE_i$ where
$\theta_1,\ldots,\theta_m$ are the eigenvalues of $A$ and
$E_1,\ldots,E_m$ are the orthogonal projections on the corresponding
eigenspaces. We may assume that $\theta_1=k$ so that
$E_1=\frac{1}{v}J_v$.
Then $A-aI=\sum_{i=1}^m(\theta_i-a)E_i$,
$A^\text{t}-aI=\sum_{i=1}^m(\overline{\theta_i}-a)E_i$ and so  
$(A-aI)(A^\text{t}-aI)=
\sum_{i=1}^m(\theta_i-a)(\overline{\theta_i}-a)E_i=
(k-a)^2E_1+b^2\sum_{i=2}^mE_i=(k-a)^2\frac{1}{v}J+b^2(I-\frac{1}{v}J).$
This equation is equivalent to 
$AA^\text{t}=(b^2-a^2+\frac{(k-a)^2-b^2}{v})I+(a+\frac{(k-a)^2-b^2}{v})(A+A^\text{t})+
\frac{(k-a)^2-b^2}{v}(J-I-A-A^\text{t})$. If $G$ is not undirected and
not a tournament then
clearly, $\lambda$ and $\mu$ are integers and then the theorem follows from
Proposition~\ref{prop:matdef}. If $G$ is a tournament then $\lambda$
is an integer, there are infinitely many choices for $(a,b)$, and
$\mu$ is arbitrary. 
If $G$ is undirected then it is strongly regular.
\qed

\begin{prop} \label{prop:eigenvals}
Suppose that $A$ is the adjacency matrix of a connected
$\nrd(v,k,\lambda,\mu)$.
Then
\begin{enumerate}
\item $k$ is an eigenvalue of multiplicity 1.
\item If $\theta$ is an eigenvalue of $A$ then ${\overline{\theta}}$ is
  an eigenvalue of the same multiplicity.
\item The spectrum of $A$ is completely determined by the spectrum of
  $A+A^\text{t}$ and the parameters $(v,k,\lambda,\mu)$.
\item if the digraph is not an undirected strongly regular graph then
  $A$ has at least one non-real eigenvalue. 
\end{enumerate}
\end{prop}

\proof
1. This is true for any connected $k$-regular digraph, see~\cite{GodRoy}.\\
2. follows from the introduction to this section.\\
3. Suppose that $\tau$ is an eigenvalue of $A+A^\text{t}$ of
multiplicity $m$. Since $A+A^\text{t}$ is symmetric, $\tau\in
\mathbb{R}$. If $|\frac{\tau}{2}-(\lambda-\mu)|=\sqrt{\eta}$ then
$\frac{\tau}{2}$  is an
eigenvalue of $A$ of multiplicity $m$. Otherwise there are exactly two
numbers $\theta$ and ${\overline{\theta}}$ on the circle with centre
$\lambda-\mu$ and radius $\sqrt{\eta}$ so that
$\theta+{\overline{\theta}}=\tau$. Then $\theta$ and ${\overline{\theta}}$
are eigenvalues of $A$ of multiplicity $\frac{m}{2}$.\\
4. If all eigenvalues of $A$ are real then since $A$ is normal
it follows that $A$ is selvadjoint and thus symmetric. But $A$ has
directed edges.\qed

\bigskip
\noindent {\bf Remark.} If $A$ is the adjacency matrix of a digraph
$\Gamma$ without undirected edges then $A+A^\text{t}$ is the adjacency
matrix of the underlying undirected graph of $\Gamma$, i.e., the graph
obtained by replacing each directed edge by an undirected edge.

\bigskip
This seems to be all that we can say in general about the the spectrum
of a normally regular digraph. But for $\mu=0$ we can at least describe the spectrum for
the most important class of normally regular digraphs.

\begin{thm}
Suppose that $T$ is a doubly regular tournament so that
$G=\mathcal{D}(T)$ is an $\nrd(v,k,\lambda,0)$. Then the eigenvalues of
$\mathcal{K}_s(G)$ are 
$$sk+s-1,\quad -1,\quad s-1+is\sqrt{k},\quad s-1-is\sqrt{k}$$ 
with multiplicities 
$$1,\quad sv-k-2,\quad \lambda+1,\quad \lambda+1.$$
\end{thm}

\proof
First we consider the eigenvalues of the adjacency matrix $A$ of $G$.
Then $A+A^\text{t}$ is the adjacency matrix of an imprimitive strongly
regular graph. This graph has eigenvalues $2k$, $0$ and $-2$ with
multiplicities $1$, $k+1$ and $k$. Thus if $\theta\neq k$ is an
eigenvalue of $A$ then $\theta+\overline{\theta}\in\{0,-2\}$. By
equation~\ref{eq:eigenvals}, $|\theta-\lambda|=\lambda+1$, as
$\eta=k-\mu+(\mu-\lambda)^2=2\lambda+1+\lambda^2$. If
$\theta+\overline{\theta}=-2$ then $\theta=-1$. The multiplicity is
$k$. If $\theta+\overline{\theta}=0$ then $\theta =\pm i
\sqrt{k}$. These two eigenvalues have multiplicity
$\frac{1}{2}(k+1)=\lambda+1$. 

Thus $A+I$ has eigenvalues $k+1$, $0$ and $1\pm i\sqrt{k}$. The matrix
$(A+I)\otimes J_s-I_{vs}$ is an adjacency matrix of
$\mathcal{K}_s(G)$. We first compute the eigenvalues of $(A+I)\otimes
J_s$ and then subtract 1.

For each eigenvector $x\in \mathbb{C}^v$ of $A+I$ with eigenvalue
$\theta$ we can replace each entry $x_i$ with $s$ entries equal to
$x_i$ to get an eigenvector in $ \mathbb{C}^{vs}$ of $(A+I)\otimes
J_s$ with eigenvalue $s\theta$. Furthermore we can get $v(s-1)$
orthogonal eigenvectors with eigenvalue $0$, by taking one of the
blocks to be orthogonal to $(1,\ldots,1)^\text{t}\in \mathbb{C}^s$ and
all other entries 0. \qed

\medskip
Note that the two non-isomorphic $\nrd(16,7,3,0)$ mentioned in
Example~\ref{Q16} have different spectra.

\section{Relation to association schemes}

An asymmetric normally regular digraph may have the additional
property that $A^2$ (where $A$ is the adjacency matrix) can be
expressed as linear combination of $A$, $A^\text{t}$, $I$ and $J$. In that case the
digraph is related to an association scheme.

\begin{defi} \label{defasssch}
Let $X$ be finite set and let $\{R_0,R_1,\ldots,R_d\}$ be a partition
of $X\times X$. Then $\mathcal{X}=(X,\{R_0,R_1,\ldots,R_d\})$
is an association scheme with $d$ classes if the following conditions
are satisfied
\begin{itemize}
\item $R_0=\{(x,x)\mid x\in X\}$,
\item for each $i\in \{0,\ldots,d\}$ there exists $i'\in
  \{0,\ldots,d\}$ such that $$R_{i'}=\{(x,y)\mid (y,x)\in R_i\},$$
\item for each triple $(i,j,k)$, $i,j,k\in\{0,\ldots,d\}$ there exist
  a number $p_{ij}^k$ such that for all $x,y\in X$ with $(x,y)\in R_k$
  there are exactly $p_{ij}^k$ elements $z\in X$ so that $(x,z)\in
  R_i$ and $(z,y)\in R_j$.
\end{itemize}
If $i=i'$ for all $i$ then the association scheme is called symmetric,
otherwise it is non-symmetric.
\end{defi}

The relation $R_i$, $i=1,\ldots,d$ can be considered as an undirected
graph if $i=i'$ and as a directed graph if $i\neq i'$. 

It is well-known that an undirected graph is strongly regular if and
only if it is a relation of a symmetric association scheme with two
classes.
And a directed graph is a doubly regular tournament if and only if it
is a relation of a non-symmetric association scheme with two
classes.

For a general introduction to association schemes we refer to Bannai and
Ito~\cite{BanI84}. Goldbach and Classen~\cite{GoldC-DM} have studied non-symmetric
association schemes with three classes and in~\cite{GolC96} they
describe the structure of non-symmetric association
schemes with three classes that are imprimitive, i.e. at least one the
graphs $R_1,R_2,R_3$ is disconnected. For tables of feasible parameter
sets see~\cite{J-3-scheme}. 

\begin{prop}

\begin{enumerate}
\item If $(X,\{R_0,R_1,R_2,R_3\})$ is an association scheme with
  $1'=2$ then $R_1$ is an asymmetric normally regular digraph.
\item If $(X,\{R_0,R_1,R_2,R_3,R_4\})$ is an association scheme with
  $1'=2$ and $3'=4$ then $R_1$ is an asymmetric normally regular digraph.
\end{enumerate}
\end{prop}

\proof
We prove case 2. Case 1 is similar. The graph $R_1$
is regular with degree $p_{12}^0$. Suppose that
$x$ and $y$ are adjacent in $R_1$. We may assume that $(x,y)\in
R_1$. Then the number of common out-neighbours of $x$ and $y$ is $p_{12}^1$.
Suppose now that $x$ and $y$ are non-adjacent. We may assume that
$(x,y)\in R_3$, since otherwise $(x,y)\in R_4$ and then $(y,x)\in R_3$.
Then the number of common out-neighbours of $x$ and $y$ is $p_{12}^3$.
\qed

\bigskip
It follows from Proposition~\ref{prop:eigenvals} that the adjacency
matrix of a normally regular digraph has at least three distinct eigenvalues. 
The normally regular digraphs constructed from non-symmetric
association schemes with 2, 3 or 4 classes have 3, 4 and 5 distinct
eigenvalues, respectively.

We now consider normally regular digraphs where the number of distinct
eigenvalues is either 3, 4 or 5, and try to construct association schemes.

In the following proofs it is easier to work with a reformulation of
the definition of association schemes in terms of matrices.

\begin{prop} \label{mat-scheme}
Let $R_0,\ldots,R_d$ be relations on a set $X$, with adjacency
matrices $A_0,\ldots, A_d$. Let $\mathcal{A}$ be the vector space
spanned by $\{A_0,\ldots, A_d\}$.
Then $(X,\{R_0,\ldots,R_d\})$ is an
association scheme if and only if $A_0+\ldots +A_d=J$ and
\begin{itemize}
\item $I\in \{A_0,\ldots, A_d\}$, (say $I=A_0$),
\item $A_i^\text{t}\in \{A_0,\ldots, A_d\}$, for all $i$, and 
\item $\mathcal{A}$ is closed under matrix multiplication.
\end{itemize}
\end{prop}

In fact, $A_iA_j=\sum_k p_{ij}^k A_k$.

\begin{thm}
Suppose that $\Gamma$ is a connected $\nrd(v,k,\lambda,\mu)$ with exactly
three distinct eigenvalues.
Then either 
\begin{itemize}
\item $\Gamma$ is an undirected strongly regular graph or
\item $\Gamma$ is a doubly regular tournament and the eigenvalues are
$k$, $-\frac{1}{2}+i\sqrt{\lambda+\frac{3}{4}}$ and $-\frac{1}{2}-
i\sqrt{\lambda+\frac{3}{4}}$, with multiplicities $1$,
$k$ and $k$ respectively.
\end{itemize}
In both cases we have an association scheme with 2 classes.
\end{thm}

\proof
Let $k$, $\theta$ and $\tau$ be the eigenvalues and let $E_k$,
$E_\theta$ and $E_\tau$ be the orthogonal projections on the
eigenspaces. Then
$$I=E_k+E_\theta+E_\tau,$$
$$ J=v E_k$$
$$A=kE_k+\theta E_\theta+\tau E_\tau,$$
and the adjoint of $A$
is
$$A^\text{t}=kE_k+\overline{\theta} E_\theta+\overline{\tau} E_\tau.$$
It follows that $I,J, A$ and $A^\text{t}$ are linearly dependent, so
there exists rational numbers $a,b,c,d$ not all 0 so that
$$aA+bA^\text{t}=c(J-I)+dI.$$
Clearly $d=0$. If $c=0$ then $A=A^\text{t}$. Otherwise either
$A+A^\text{t}=J-I$ or $A=A^\text{t}=J-I$.
Thus $I, A, A^\text{t}$ satisfy the properties required in
Proposition~\ref{mat-scheme}. 

The spectrum can be computed as described in
Proposition~\ref{prop:eigenvals}. \qed

\begin{thm} \label{thm:4eigenvals}
Suppose that $\Gamma$ is a connected asymmetric
$\nrd(v,k,\lambda,\mu)$ with exactly four distinct eigenvalues.

Then $\Gamma$ is a relation of a non-symmetric association scheme
with three classes.
\end{thm}

\proof
Let $A$ be the adjacency matrix of $\Gamma$. 
Let
$\mathcal{A}=\text{span}\{I,A,A^\text{t},J-I-A-A^\text{t}\}$.
We need to show that $\mathcal{A}$ is closed under
matrix multiplication. Let $k,\tau,\theta$ and ${\overline{\theta}} $ be the
eigenvalues of $A$. Let $E_k, E_\tau, E_\theta$ and
$E_{\overline{\theta}}$, respectively, be the orthogonal projections on
the corresponding eigenspaces.
Then
$$I=E_k+E_\tau+ E_\theta+E_{\overline{\theta}},$$
$$J=vE_k,$$
$$A=kE_k+\tau E_\tau+\theta
E_\theta+{\overline{\theta}}E_{\overline{\theta}},$$
and 
$$A^\text{t}=kE_k+\tau E_\tau+{\overline{\theta}}
E_\theta+\theta E_{\overline{\theta}}.$$
Thus $\{E_k,E_\tau,E_\theta,E_{\overline{\theta}}\}$ is a basis of
$\mathcal{A}$. Since these projections are idempotents and the product
of distinct projections is 0, $\mathcal{A}$ is closed under
multiplication. \qed

\bigskip
\noindent {\bf Remark.} We did not use the fact that $\Gamma$ is a
normally regular digraph in the proof of Theorem~\ref{thm:4eigenvals}.
In fact we proved that if $\Gamma$ is a regular connected graph
without undirected edges whose adjacency matrix is normal and has
exactly four distinct eigenvalues then $\Gamma$ is a relation of a
non-symmetric association scheme with three classes.

If an asymmetric normally regular digraph has exactly five distinct eigenvalues
then it may have either three real eigenvalues and one pair of
complex conjugate eigenvalues or else it has one real eigenvalues and
two pairs complex conjugate eigenvalues. In the latter cases it seems
likely that the graph is a relation of a non-symmetric association
scheme with four classes. We can only prove the following.

\begin{prop}
Suppose that $\Gamma$ is a connected asymmetric
$\nrd(v,k,\lambda,\mu)$ with exactly five distinct eigenvalues
$k,\theta,{\overline{\theta}},\tau,\overline{\tau}$.

Then $\Gamma$ is an orientation of a strongly regular graph.
\end{prop}

\proof
The underlying undirected graph is regular and has exactly three
distinct eigenvalues $2k$, $\theta+{\overline{\theta}}$ and
$\tau+\overline{\tau}$. Thus it is strongly regular, see~\cite{GodRoy}. \qed

\medskip

Conversely, it follows from Proposition~\ref{prop:eigenvals} that if a
normally regular digraph is an orientation of a strongly regular graph
then the number of distinct eigenvalues is either four or five.

\section{Group divisible partitions\label{groupdivisible}}

We start with the definition of two types of partitions of the vertex set.

Suppose that the vertex set of a normally regular digraph is
partitioned in sets $V_1,\ldots ,V_m$.
Then we say that $V_1,\ldots ,V_m$ is an
equitable partition if there exists constants
$c_{ij}, d_{ij}$ for $i,j\in\{1,\ldots,m\}$ so that for every vertex $x\in
V_i$, $|x^+\cap V_j|=c_{ij}$ and for every vertex $y\in V_j$,
$|y^-\cap V_i|=d_{ij}$. If $|V_i|=|V_j|$ then $c_{ij}=d_{ij}$. 
We say that $C=(c_{ij})_{i,j=1,\ldots,m}$ is the quotient matrix of
the equitable partition, see~\cite{GodRoy}.

Let $G$ be an asymmetric $\nrd(v,k,\lambda,\mu)$. Then we say that $G$ is group
divisible if $G$ is a multipartite tournament, i.e., if $V(G)$ can be
partitioned in sets $V_1,\ldots,V_r$ such that there is an edge
between $x\in V_i$ and $y\in V_j$ if and only if $i\neq j$.

Since $G$ is regular the sets $V_i$ all have the same size, say
$|V_i|=s=v-2k$, for $i=1,\ldots,r$. Then $v=rs$ and
$k=\frac{1}{2}(r-1)s$. We assume that $s>1$.

The adjacency matrix of a group divisible normally regular digraph with $\mu\neq \lambda$ is also the incidence
matrix of a group divisible design, see~\cite{Beth}

\newpage
\begin{lemma}
Let $G$ be an asymmetric $\nrd(v,k,\lambda,\mu)$ with a partition
$V_1,\ldots,V_r$ of the vertex set so that two vertices are adjacent if
and only if they are in different cells.

Then $V_1,\ldots,V_r$ is an equitable partition.
\end{lemma}

\proof
For $x\in V_i$ let $c_{ij}(x)= |x^+\cap V_j|$. We need to prove that
$c_{ij}(x)$ does not depend on $x$.
Let $x\in V_i$ and $y\in V_j$.
We count the vertices in $S=\{z\mid x\pil z\pil y\}$ in two ways.
The number of out-neighbours of $x$ outside $V_j$ is
$k-c_{ij}(x)$. $\lambda$ of these out-neighbours are common
out-neighbours of $x$ and $y$. The remaining $k-c_{ij}(x)-\lambda$
vertices are in $S$. Similarly, $y$ has $k-(|V_i|-c_{ji}(y))=k-(s-c_{ji}(y))$
in-neighbours outside $V_i$. $k-(s-c_{ji}(y))-\lambda$ vertices are
in $S$. Thus $s=c_{ij}(x)+c_{ji}(y)$, for all $x\in V_i$. \qed

\bigskip

\begin{prop}
Let $A$ be the adjacency matrix of a group divisible normally regular
digraph. Then
\begin{itemize}
\item $A$ has either 4 or 5 distinct eigenvalues.
\item If $A$ has 4 distinct eigenvalues then the graph is a relation of a
  non-symmetric imprimitive association scheme with three classes.
\item If $A$ has 5 distinct eigenvalues then $r$ and $s$ are odd.
\end{itemize}
\end{prop}

\proof
Since $A+A^\text{t}$ is a strongly regular graph it has three distinct
eigenvalues. Then by Theorem~\ref{thm:eigenval} and
Proposition~\ref{prop:eigenvals}, $A$ has either 4 or 5 eigenvalues.

It follows from Theorem~\ref{thm:4eigenvals} that if $A$ has 4
eigenvalues then the graph is a relation of a non-symmetric
association scheme with three classes. 

It follows from  Proposition~\ref{prop:eigenvals} that if $A$ has 5
distinct eigenvalues then eigenvalue $k$ has multiplicity 1 and the
other 4 eigenvalues have pairwise the same multiplicity. Thus the
number of vertices is odd, and so $r$ and $s$ are odd. \qed

\bigskip

I conjecture that if a group divisible normally regular digraph with an
odd number of vertices exists then it satisfies $\mu=k$ and then by
Theorem~\ref{myk} it is a relation of a non-symmetric imprimitive
association scheme with three classes. 
\begin{conj}
Any group divisible normally regular digraph is a relation of a
non-symmetric imprimitive association scheme with three classes. 
\end{conj}

\begin{ex}
The parameters of a group divisible normally regular
digraph must satisfy that $v-2k$ divides $v$. This is satisfied by
$(v,k,\lambda,\mu)=(16,6,2,2)$. There are four asymmetric normally
regular digraphs with these parameters. Two of these are group
divisible and thus are relations of an association scheme.

One of these is a Cayley graph 
$$\cay(\mathbb{Z}_4\times \mathbb{Z}_4,
\{(0,3),(1,3),(2,1),(3,0),(3,2),(3,3)\}).$$ 
The independent sets of
vertices in this digraph are the cosets of the subgroup
$\{(0,0),(0,2)$, $(2,0),(2,2)\}$.

One of the normally regular digraphs with these parameters that is not
group divisible has vertex set $\{a_i,b_i\mid i\in \mathbb{Z}_8\}$ and
edges 
$$a_i \pil a_{i+1}, a_{i+2}, b_i, b_{i+1}, b_{i+4}, b_{i+6}, \ \
i\in\mathbb{Z}_8, $$
$$b_i \pil b_{i-1},b_{i-2}, a_{i-2},a_{i-3}, a_{i-5}, a_{i-7}, \ \
i\in\mathbb{Z}_8.$$ 

Thus group divisibility is not determined by the parameters.
\end{ex}

\section{Combinatorial results for small $\lambda$}

In this section we use combinatorial methods to prove non-existence
for certain parameter sets where $\lambda$ is small.
If $\lambda$ is so small that $2\lambda-\mu$ is negative then only the
asymmetric case need to be considered.

\begin{thm} \label{lalille}
If there exist an asymmetric normally regular digraph with parameters $(v,k,\lambda,\mu)$ where
$2\mu>k+\lambda$ then the graph is group-divisible and $v-2k$ divides
$v$.
\end{thm}
\proof
Suppose that $G$ is an $\nrd(v,k,\lambda,\mu)$. Let $x$ be a vertex in
$G$. Let $y$ and $z$ be vertices in $V(G)-\{x\}-x^+-x^-$. Then
$x^+\cap y^+$ and $x^+\cap z^+$ each consist of $\mu$ vertices in the
set $x^+$ of $k$ vertices. Thus $|y^+\cap z^+|\geq 2\mu-k>\lambda$ and
so $y$ and $z$ are nonadjacent. It follows that every vertex in $G$
belongs to a unique independent set of $v-2k$ vertices and so  $G$ is
group divisible with $\frac{v}{v-2k}$ groups. \qed

\begin{cor}
Suppose that $2\lambda-\mu<0$, $2\mu>k+\lambda$ and $v-2k$ does not divide
$v$ then an $\nrd(v,k,\lambda,\mu)$ does not exist.
\end{cor}

\begin{thm} \label{comb}
If a normally regular digraph with $\lambda=0$, $\mu\neq k$ and $\mu\geq 2$ exists then 
$k\geq 2 \mu+\frac{1}{2}+\sqrt{2\mu+\frac{1}{4}}$.
\end{thm}

From Theorem~\ref{myk} we know that an $\nrd(v,k,0,k)$ is obtained from a
directed triangle by replacing each vertex by $k$ vertices.

\medskip
\noindent {\bf Proof\ }
Suppose that $G$ is an $\nrd(v,k,0,\mu)$ with $\mu\neq k$.
Since $2\lambda-\mu<0$ any such normally regular digraph is asymmetric.
In this proof we use notation $U_x=V(G)-\{x\}-x^+-x^-$ for a vertex $x$ in
$G$. By equation~\ref{param-eq1} we have
\begin{equation} \label{eq:la0}
|U_x|=\frac{k(k-1)}{\mu}
\end{equation}

By Proposition~\ref{mu-interval}, $k>\mu$.

\noindent {\bf Claim 1:} $k> 2\mu$.

{\bf Proof\ } Let $x$ be any vertex in $G$.
By Theorem~\ref{lalille}, $U_x$ is an independent set if $\mu<k<2\mu$.

So suppose that $k=2\mu$ and $k>2$. By equation~\ref{eq:la0}, $|U_x|=2k-2$.
Suppose that $y,z\in U_x$
and $y$ dominates $z$. $x^+\cap y^+$ and $x^+\cap z^+$ are disjoint
sets (as $\lambda=0$) of cardinality $\mu$. Thus their union is $x^+$.
$z$ has $\mu$ in-neighbours in $x^-$, no in-neighbours in $x^+$ and
thus $\mu$ in-neighbours in $U_x$. Let $y'\in U_x$ be another vertex
dominating $z$. Then $y'$ and $z$ have no common out-neighbours in
$x^+$, i.e. $x^+\cap (y')^+= x^+\cap y^+$. Thus $y$ and $y'$ have at
least $\mu+1$ common out-neighbours, a contradiction.

Thus $U_x$ is an independent set.

Let $z\in x^+$. Every vertex other than $x$ dominating $z$ belong to
$U_x$. As $U_x$ is independent, no vertex dominates both $z$ and $y\in U_x$
and so $z$ is adjacent to every vertex in $U_x$ (but to no other vertex
in $x^+$).
Thus $z$ is adjacent to $2k-1-|U_x|$ vertices in $x^-$.
By equation~\ref{eq:la0} and $\mu<k\leq 2\mu$, we have
$k>2k-1-|U_x|\geq 1$, and so there is a vertex $w\in x^-$ adjacent to $z$ and a
vertex $u\in x^-$ not adjacent to $z$. Then $w$ and $z$ are adjacent
vertices in $U_u$, a contradiction.

\noindent{\bf Claim 2:} $k>2\mu+1$.

{\bf Proof\ } Suppose that $k=2\mu+1$ and let $x$ be a vertex in $G$. By
equation~\ref{eq:la0}, $|U_x|=2k$. Since $G$ is regular and
$|U_x\cup\{x\}|>|x^+\cup x^-|$, $U_x$ cannot be an independent set.
Let $y,z\in U_x$ so that $y\rightarrow z$. $y$ and $z$ have no common
out-neighbours in $x^+$ so there a unique vertex $w\in x^+$ which is
not dominated by $y$ or $z$. The in-neighbours of $z$ are $\mu$
vertices in $x^-$, possibly $w$, and at least $\mu$ vertices in $U_x$.
Let $y'\neq y$ be a vertex in $U_x$ dominating $z$. Since $y'$ and $z$
have no common out-neighbours in $x^+$ and $y$ and $y'$ have only
$\mu$ common out-neighbours, $y'$ dominates $w$. Since $\lambda=0$,
$w$ does not dominate $z$, and so $z$ has $\mu+1$ in-neighbours in
$U_x$.
We have now shown that in the
graph spanned by $U_x$ any vertex has in-degree either 0 or $\mu+1$
and, by symmetry, it has out-degree either 0 or $\mu+1$. Thus the
in-neighbours of $z$ in $U_x$ has out-degree $\mu+1$.
Any two in-neighbours of $z$ in $U_x$ have at least $\mu-1$
common out-neighbours in $(x^+\cap y^+)\cup\{w\}$. Thus $z$ is their
only common out-neighbour in $U_x$.
Counting the vertices in $U_x$ we have
$$4\mu+2=|U_x|\geq 1+(\mu+1)+(\mu+1)\mu,$$
and so $\mu=2$.

Now $|U_x|=10$ and we have at least 7 vertices of in-degree 3 in $U_x$
and, by symmetry, at least 7 vertices of out-degree 3. So there is a
vertex with out-degree and in-degree 3. We may assume that $z$ is such
a vertex. Then $z$ dominates a vertex which is also dominated by an
in-neighbour of $z$. This is a contradiction to $\lambda=0$.
This proves claim 2.

\bigskip

Let $r=k-2\mu$. Then $r\geq 2$, by Claim 2. By equation~\ref{eq:la0}, $\mu$
divides $k(k-1)=(2\mu+r)(2\mu+r-1)=\mu(4\mu+4r-2)+r(r-1)$. Thus $\mu$
divides $r(r-1)$ and so $r^2-r=s\mu$, for some positive integer $s$.
Then $r=\frac{1}{2}+\sqrt{s\mu+\frac{1}{4}}$.

If $s=1$ then $\mu=r(r-1)$ and $k=2\mu+r=r(2r-1)$. From
equation~\ref{eq:la0},
we see that $v=1+2k+\frac{k(k-1)}{\mu}=2k+4r^2$ is even
and so $\eta$ is a square,
by Theorem~\ref{BRC}. But $\eta=k-\mu+\mu^2=r^2((r-1)^2+1)$ cannot be a square.
Thus $s\geq 2$. This proves the theorem.\hfill $\square$

\section{Normally regular digraphs as quotient graphs}

\subsection{Subplane partition}
\label{quotient}

Fossorier, Je\v{z}ek, Nation and Pogel~\cite{ordinary} considered
partition of a projective plane of order $n$ into subplanes
$\pi_1,\ldots,\pi_v$ of order $q$, $v=\frac{n^2+n+1}{q^2+q+1}$. They
say that such a partition is 
ordinary if for each pair $(i,j)$ either each point of $\pi_i$ is
incident with a line of $\pi_j$ or no point of $\pi_i$ is
incident with a line of $\pi_j$. 

For an ordinary partition of a projective plane they consider the
quotient graph with vertices $\pi_1,\ldots,\pi_v$ and an edge
$\pi_i\pil\pi_j$ if the points of $\pi_i$ are incident with lines of
$\pi_j$. They proved that this quotient graph is what they called an
ordinary graph. This is a normally regular digraph in our terminology.

\begin{thm} [Fossorier, Je\v{z}ek, Nation and Pogel~\cite{ordinary}]
\label{subplanepart}
If a projective plane of order $n$ has an ordinary partition into
projective planes of order $q$ then the quotient graph is a normally
regular digraph with $(v,k,\lambda,\mu)=(\frac{n^2+n+1}{q^2+q+1},n-q,
q^2,q^2+q+1)$. 
\end{thm}

For a partition into Baer subplanes, i.e., $n=q^2$, the quotient
graph is a complete undirected graph.

Theorem~\ref{subplanepart-desargues} below describes the special case
of this theorem where we consider desarguesian planes. 
Theorem~\ref{k13} may also be seen as a special case of
Theorem~\ref{subplanepart} where $q=1$. 

\subsection{Bipartite graphs of diameter 3}

Delorme, J{\o}rgensen, Miller and Pineda-Villavicencio~\cite{diam3}
considered a similar quotient graph construction. 
In this paper we considered bipartite $q+1$ regular graphs with
diameter 3 and with $2(q^2+q)$ vertices.
(The largest possible bipartite $q+1$ regular graph with diameter 3
has $2(q^2+q+1)$ vertices and it appears only as incidence graph of a
projective plane of order $q$.)
In such graphs the vertices are partitioned into cycles of length 4.
It is proved that the graph obtained by directing all edges from one
bipartition class to the other and then identifying each 4-cycle to a
vertex is a normally regular digraph with
$(v,k,\lambda,\mu)=(\frac{q^2+q}{2},q-1,0,2)$.

This was our original motivation for studying normally regular digraphs. 

\section{Constructions}

In this section we give a number constructions of families of normally
regular digraphs. Most of these constructions use Cayley graphs of
abelian groups. 

\subsection{Asymmetric Cayley graph constructions}

The first construction uses a partition of a projective plane into
triangles. If a triangle is considered to be a ``subplane'' of order 1
then this is a special case of the construction in
Theorem~\ref{subplanepart}. 

This construction was also found by de Resmini and
Jungnickel~\cite{ResJun} as an example of what they call a failed
symmetric design.

\begin{thm} \label{k13}
Let $k$ be a multiple of 3 such that $k+1$ is a prime power. Then
there exists $S\subset \mathbb{Z}_v$, $v= \frac{k^2+3k+3}{3}$ so that 
$\cay(\mathbb{Z}_v,S)$ is an asymmetric normally regular digraph with 
parameters $(v,k,1,3)$.
\end{thm}

\noindent
{\bf Proof\ } When $\lambda=1$ and $\mu=3$, $\eta=k+1$, which is assumed
to be a prime power. By equation~\ref{eq:veq2},
$v=\frac{\eta^2+\eta+1}{3}$.

By Singer's theorem~\cite{Singer} there exist a cyclic planar
difference set of order $\eta$, i.~e. a subset $D$ of $\mathbb{Z}_{3v}$ with
$|D|=\eta+1$ such that each non-zero element of $\mathbb{Z}_{3v}$ is a
difference of exactly one ordered pair of elements in $D$. In particular there
is a unique pair of difference $v$. By adding a constant to $D$ if
necessary, we may assume that $v,2v\in D$. Let
$D'=D\setminus\{v,2v\}$, and let $S\subseteq \mathbb{Z}_v$ be the numbers
congruent to numbers in $D'$ modulo $v$.

As $v$ is not a difference in $D'$, $|S|=\eta-1=k$.

Suppose that $x$ and $-x$ are both in $S$. Then for some
$a,b\in\{0,1,2\}$, $x+av,-x+bv\in D'$. Choose $i,j\in\{1,2\}$ so that
$a+b\equiv i+j$ (mod 3). Then we have two equal differences in $D$
$$(x+av)-iv=jv-(-x+bv)$$
a contradiction. Thus $S\cap -S=\emptyset$.

If $x\in \mathbb{Z}_v$ is congruent (mod $v$) to a difference of elements
in $D$, one of which is $v$ or $2v$, then either $x\equiv a-iv$ or $x\equiv
iv-a$ for $a\in D$ and $i\in \{1,2\}$, i.~e. $x$ or $-x\in S$.
Conversely if $x\in S$ or $-x\in S$ then $x$ is in exactly to ways
congruent (mod $v$) to
a difference of two elements in $D$ one of which is $v$ or
$2v$.

Let $x\in \mathbb{Z}_v\setminus(S\cup -S\cup\{0\})$. Then each of $x$,
$x+v$ and $x+2v$ can be written in exactly one way as a difference of
elements in $D$, in fact in $D'$. Thus $x$ can be written in exactly
three ways as a difference of elements in $S$.

If $x\in S\cup -S$ then only one of the three pairs of elements in $D$
whose difference is congruent to $x$ (mod $v$) is in $D'$.
Thus $x$ can be written in exactly one way as a difference of elements
in $S$.

Hence $\cay(\mathbb{Z}_v,S)$ is an $\nrd(v,k,1,3)$. \qed

\bigskip

In the next theorem we construct a family of Cayley graphs of abelian
but not necessarily cyclic groups. It is well-known that this
digraph is one of the classes of a (so-called cyclotomic) association
scheme with four classes. 

\begin{thm} \label{quartic}
Suppose that $v$ is a prime power, $v\equiv 5 \mod 8$.
Let $D$ denote the following subset of GF[$v$]:
$$ D=\{x^4\mid x\neq 0\}$$
Then the Cayley graph of the additive group of GF[$v$] generated by
$D$ is a normally regular digraph with $v=4k+1=8(\mu+\lambda)+5$.
\end{thm}

\noindent {\bf Proof\ }  As
$v\equiv 5 \mod 8$, the set $D$ has cardinality $k=\frac{v-1}{4}$
and $-1\not\in D$. Thus
$D\cup -D$ is the set of squares in $GF[v]$. This means that the cosets
of the subgroup (of the multiplicative group) $D$ are $D$, $-D$, $R$,
and $-R$ for some set $R$. Let $D=\{1,q_2,\ldots,q_k\}$.
Then $qq_i-q$, $2\leq i\leq k$, $q\in D$ is the set of differences, we
want to consider. For a fixed $i$, every element in the coset to which
$q_i-1$ belongs appears exactly once as a difference $qq_i-q$, $q\in
D$. This means that if among the differences $q_2-1,\ldots,q_k-1$,
the number of elements in $D,-D,R$ and $-R$,
are $\lambda_1,\lambda_2,\mu_1$ and $\mu_2$,
respectively, then among all differences of distinct element of $D$ an
element appears $\lambda_1,\lambda_2,\mu_1,$ or $\mu_2$ times according to
whether it belongs to $D,-D,R,$ or $-R$. Since for every $x$, $x$ and $-x$
appears as a difference the same number of times,
$\lambda_1=\lambda_2$ and $\mu_1=\mu_2$.
\qed

\bigskip

The only known infinite family of primitive non-symmetric association
schemes with three classes is a family constructed by Liebler and
Mena~\cite{LM}. 
For every $s=2^n$, they constructed a so-called distance regular
digraph of girth 4 and degree $s(2s^2-1)$, as a Cayley digraph
of ${\mathbb{Z}}_4 \times \ldots \times{\mathbb{Z}}_4$.
Their graph is in fact an $\nrd(4s^4,s(2s^2-1),2s(s-1),s(s-1))$.

\bigskip
Some of the normally regular digraphs constructed in the next two
subsections are also asymmetric.

\subsection{Construction from desarguesian planes}

We will now consider the subplane partition described in
Section~\ref{quotient} for desarguesian projective planes.

\begin{thm} \label{subplanepart-desargues}
Let $q$ be a prime power and let $r\geq 2$ be an integer not divisible
by 3. Let $v=\frac{q^{2r}+q^r+1}{q^2+q+1}$. Then there exists a set
$S\subset \mathbb{Z}_v$ so that $\cay(\mathbb{Z}_v,S)$ is a normally
regular digraph with parameters $(v,q^r-q,q^2,q^2+q+1)$.
\end{thm}

\proof
Let GF[$q^{3r}$] be the field with $q^{3r}$ elements and with
primitive element $\alpha$. Then GF[$q^{r}$] and GF[$q^{3}$]
are subfields and their intersection is GF[$q$] as 3 does not divide
$r$.

Let $\beta=\alpha^{(q^{3r}-1)/(q^3-1)}$. Then $\beta$ is a primitive
element of GF[$q^{3}$] and $\beta\notin $  GF[$q^r$]. Thus when
GF[$q^{3r}$] is considered as a 3 dimensional vector space over
GF[$q^{r}$] then vectors $1$ and $\beta$ span a 2 dimensional subspace
$U$.
Let $u=\frac{q^{3r}-1}{q^r-1}=q^{2r}+q^r+1$ and let
$D=\{i\in\mathbb{Z}_u\mid \alpha^i\in U\}$. Then by Singer's
theorem~\cite{Singer}, $D$ is a planar difference set in
$\mathbb{Z}_u$.

Similarly, we may consider  GF[$q^{3}$] as a 3 dimensional vector
space over GF[$q$]. In this space the vectors $1$ and $\beta$ span a 2
dimensional subspace $W$.
Let $w=\frac{q^{3}-1}{q-1}=q^{2}+q+1$ and let
$T=\{i\in\mathbb{Z}_w\mid \beta^i\in W\}$. Again $T$ is a planar
difference set in $\mathbb{Z}_w$. As $\frac{q^r-1}{q-1}$ and $w$ are
coprime, multiplication by $\frac{q^r-1}{q-1}$ is an automorphism of
$\mathbb{Z}_w$ and so $T'=\{\frac{q^r-1}{q-1}i \mid i \in T\}$ is
a planar difference set. Then the set $T''=\{\frac{q^r-1}{q-1}vi \mid i \in
T\}$ is a difference set in subgroup $\langle v\rangle$ of $\mathbb{Z}_u$.
This set satisfies $T''=\{i\in\mathbb{Z}_u\mid \alpha^i\in W\}$ and
$T''\subset D$, as $\beta=\alpha^{\frac{q^r-1}{q-1}v}$. 

If for some $x,y\in D$ the
difference $x-y$ is a non-zero multiple of $v$ then $x,y\in T''$. Let
$D'=D\setminus T''$. Let $S\subset \mathbb{Z}_v$ be the numbers
congruent to numbers in $D'$ modulo $v$. As multiples of $v$ are not
differences in $D'$, $|S|=q^r-q$.

Let $g\in\mathbb{Z}_v$, $g\neq 0$. Then $g$ is congruent modulo $v$ to
$q^2+q+1$ elements in $\mathbb{Z}_u$, each of which can uniquely be
written as difference $x-y$ where $x,y\in D$. If $g\in S$ then exactly
$q+1$ of these differences satisfy $y\in T''$. If $g\notin S$ then none
of the differences have $y\in T''$. Similarly, if $g\in -S$ then
exactly $q+1$ of the differences have $x\in T''$.

Thus the number of pairs $x,y\in S$ such that a nonzero element $g\in
\mathbb{Z}_v$ can be written as $g=x-y$ is $\mu=q^2+q+1$ if $g\notin
S\cup -S$, $\lambda=q^2$ if $g$ is in exactly one of the sets $S, -S$,
and $2\lambda-\mu=q^2-q-1$ if $g\in S\cap -S$.  
\qed

\bigskip
For $r=2$ the graph constructed in this theorem is a complete graph
with $v=\frac{q^4+q^2+1}{q^2+q+1}=q^2-q+1$. If $r\geq 4$ is even then
the projective plane of order $q^r$ has an ordinary partition in
subplanes of order $q^2$ and the planes of order $q^2$ have an ordinary
partition in subplanes of order $q$. Then the vertices of the normally
regular digraph are partitioned in sets of size $q^2-q+1$ spanning
complete subgraphs. Thus $S\cap -S$ contains all nonzero elements of
the subgroup of order $q^2-q+1$. 

We conjecture that these are the only elements in $S\cap -S$.

\begin{conj}
Let $S$ be as in Theorem~\ref{subplanepart-desargues}. 
\begin{itemize}
\item If $r$ is odd then $S\cap -S=\emptyset$.
\item If $r$ is even then $S\cap -S$ consists of the nonzero
  elements of the subgroup of order $q^2-q+1$. 
\end{itemize}
\end{conj}

The proof of Theorem~\ref{subplanepart-desargues}
is an algorithm for computing the set $S$. The above conjecture is
based on computations of $S$ for the following values of $(q,r)$:
(2,4), (2,5), (2,7), (2,8), (2,10), (3,4), (3,5), (4,4), (4,5), (5,4).

\medskip
\begin{ex} \label{ex:desargues}
1. For $q=2$ and $r=4$,
we get 
$$S=\{7,13,14,17,19,23,26,28,29,31,34,35,37,38\}\in \mathbb{Z}_{39}.$$
This gives a normally regular digraph with parameters $(39,14,4,7)$. 
In this particular case, the normally regular digraph is isomorphic to
the graph constructed in Corollary~\ref{cor:prod2}, with
$(s,t)=(0,3)$.\\
2. For $q=2$ and $r=5$, we get 
\begin{eqnarray*}
S &= &\{11, 17, 21, 22, 25, 29, 31, 34,42, 43, 44, 45, 49, 50, 58, 62, \\
  &  &68, 81, 84, 86, 88, 90, 91, 97, 98,100, 116, 
121, 124, 136\}\in \mathbb{Z}_{151}.
\end{eqnarray*} 
This gives an {\em
  asymmetric} normally regular digraph with parameters
$(151,30,4,7)$.\\
\end{ex}

\subsection{Product constructions}

In this section we give two constructions of normally regular digraphs
that are not asymmetric. They are products involving doubly regular
tournament and conference graphs. In some cases they are Cayley graphs.

\begin{thm}
Let $T$ be a doubly regular tournament with $4t+3$ vertices and let
$K_{2t+1}$ be the complete graph of order $2t+1$.
Then the cartesian product with vertex set $V(T)\times V(K_{2t+1})$ and
edge set $\{(x,u)\pil (y,u) \mid x\pil y \text{\ in \ } T \}
\cup\{(x,u)\dpil (x,v)\mid u\dpil v \text{\ in \ } K_{2t+1}\}$
is a normally regular digraph with parameters
$((4t+3)(2t+1),4t+1,t,1)$.
\end{thm}

\proof
Let $(x,u)$ and $(y,v)$ be vertices in the cartesian product. If $x=y$
and $u\neq v$ then  $(x,u)$ and $(y,v)$ are joined by an undirected
edge and their common out-neighbours are the remaining $2t-1$ vertices
of the form $(x,w)$.
If $x\neq y$ and $u=v$ then  $(x,u)$ and $(y,v)$ are joined by a
directed edge and their common out-neighbours are the vertices $(z,u)$
where $z$ is a common out-neighbour of $x$ and $y$ in $T$. There are
$t$ such vertices.
If $x\neq y$ and $u\neq v$ then  $(x,u)$ and $(y,v)$ are
non-adjacent. We may assume $x\pil y$ in $T$. Then $(y,u)$ is the
unique common out-neighbour of  $(x,u)$ and $(y,v)$. \qed

\medskip
If $T$ is a Paley tournament then the above construction is a Cayley
graph.

\begin{cor} \label{cor:prod1}
Let $\mathbb{F}$ be the field of $4t+3$ elements and let $Q$ be the
set of non-zero squares in $\mathbb{F}$.
Let $\mathbb{Z}_{2t+1}$ be the cyclic group of order $2t+1$.
Let $G$ be the direct product $\mathbb{F}\times \mathbb{Z}_{2t+1}$ and
let $S = \{(d,0) \mid d \in Q\} \cup \{(0,z) \mid z\neq 0\}$.

Then $\cay(G,S)$ is a normally regular digraph with
parameters $((4t+3)(2t+1),4t+1,t,1)$.
\end{cor}

\begin{thm}
Let $H$ be a conference graph with $4t+1$ vertices and let $T$ be a
doubly regular tournament with $4s+3$ vertices. Let $\Gamma$ be the
graph with vertex set $V(H)\times V(T)$ and with edge set
$$\{(u,x)\pil (v,y)\mid \text{either\ } u\dpil v \text{\ and\ } x\pil
y, \text{\ or\ } u\not\dpil v \text{\ and\ } x\bpil y\}$$
$$\cup\
\{(u,x)\dpil (u,y)\mid u\in V(H),\ x,y\in V(T)\}.$$
Then $\Gamma$ is a normally regular digraph with
parameters
$$((4t+1)(4s+3),(4t+2)(2s+1),4ts+3s+t+1,(2t+1)(2s+1)).$$
\end{thm}

\proof
Suppose that $(u,x)\pil (v,y)$ but $(u,x)\not\bpil (v,y)$
in $\Gamma$. Then either $u\dpil v$ and
$x\pil y$, or $u$ and $v$ are non-adjacent and $x\bpil y$. Suppose
that $u\dpil v$ and $x\pil y$. Let $(w,z)$ be a common out-neighbour
of $(u,x)$ and $(v,y)$. Then $w$ and $z$ satisfy one of the following
six cases.
$$w=u, y\pil z,$$
$$w=v, x\pil z, z\neq y,$$
$$u\dpil w\dpil v, x\pil z, y\pil z,$$
$$u\not\dpil w\not\dpil v, x\bpil z, y\bpil z,$$
$$u\dpil w\not\dpil v, x\pil z, y\bpil z,$$
$$u\not\dpil w\dpil v, x\bpil z, y\pil z.$$
The number of vertices $(w,z)$  in each case are $2s+1$, $2s$, $ts$,
$(t-1)s$, $ts$ and $t(s+1)$, respectively. The case where $u$ and $v$
are non-adjacent and $x\bpil y$ is similar. Thus $\lambda=4ts+3s+t+1$.

The parameters $k$ and $\mu$ are easy to compute. Now
suppose that $(u,x)\dpil (v,y)$. Then $u=v$. A common out-neighbour
$(w,z)$ of $(u,x)$ and $(u,y)$ is of one of the following three cases.
$$w=u, z\neq x,y,$$
$$w\dpil u, x\pil z, y\pil z,$$
$$w\not\dpil u, x\bpil z, y\bpil z.$$
The number of vertices of each type is $4s+1$, $2ts$ and $2ts$,
respectively.
This adds up to $4ts+4s+1=2\lambda-\mu$. \qed

\medskip
If the conference graph and the doubly regular tournament in this
theorem are both of Paley type then $\Gamma$ is a Cayley graph.

\begin{cor} \label{cor:prod2}
Let $\mathbb{F}$ and $\mathbb{E}$ be finite fields of order $4t+1$ and
$4s+3$, respectively.
Let $Q_\mathbb{F}$ and  $Q_\mathbb{E}$ be the sets of non-zero squares
in $\mathbb{F}$ and $\mathbb{E}$, respectively. Let
$R_\mathbb{F}=\mathbb{F}\setminus (Q_\mathbb{F}\cup\{0\})$ and
$R_\mathbb{E}=\mathbb{E}\setminus (Q_\mathbb{E}\cup\{0\})$. Let
$S=(Q_\mathbb{F}\times Q_\mathbb{E})\cup (R_\mathbb{F}\times
R_\mathbb{E})\cup (\{0\}\times (\mathbb{E}\setminus \{0\}))$. Then
$\cay(\mathbb{F}\times\mathbb{E},S )$ is a normally
regular digraph with parameters
$((4t+1)(4s+3),(4t+2)(2s+1),4ts+3s+t+1,(2t+1)(2s+1))$.
\end{cor}

Note that if $(4t+1)-(4s+3)=\pm 2$ then we get the difference set
with Hadamard parameters constructed by Stanton and
Sprott~\cite{StantonSprott} by adding $(0,0)$ to $S$ if $(4t+1)-(4s+3)=2$
or by talking the complement of $S$ if $(4t+1)-(4s+3)=-2$.

\end{document}